\documentclass[pdflatex,sn-apa]{sn-jnl}

\usepackage{graphicx}%
\usepackage{multirow}%
\usepackage{amsmath,amssymb,amsfonts}%
\usepackage{amsthm}%
\usepackage{mathrsfs}%
\usepackage[title]{appendix}%
\usepackage{xcolor}%
\usepackage{textcomp}%
\usepackage{manyfoot}%
\usepackage{booktabs}%
\usepackage{algorithm}%
\usepackage{algorithmicx}%
\usepackage{algpseudocode}%
\usepackage{listings}%
\usepackage{enumitem}
\usepackage{siunitx}
\theoremstyle{thmstyleone}%
\theoremstyle{thmstyletwo}%

\theoremstyle{thmstylethree}%

\begin{document}

\title{Synthetic Fluency and Epistemic Offloading in Undergraduate Mathematics in the Age of AI}


\author[1]{\fnm{Siyuan} \sur{Wang}}\email{siywang@kean.edu}
\equalcont{These authors contributed equally to this work.}

\author[1]{\fnm{Qing} \sur{Xia}}\email{qxia@kean.edu}
\equalcont{These authors contributed equally to this work.}

\author*[1]{\fnm{Qiong} \sur{Ye}}\email{yeqiong@wku.edu.cn}
\equalcont{These authors contributed equally to this work.}

\affil[1]{\orgdiv{College of Science, Mathematics and Technology}, \orgname{Wenzhou-Kean University}, \orgaddress{\city{Wenzhou}, \postcode{325060}, \state{Zhejiang}, \country{China}}}




\abstract{The rapid adoption of generative artificial intelligence (AI) tools in higher education is transforming how students engage with undergraduate mathematics, raising concerns about learning and assessment validity. This study examines the impact of AI accessibility across a two-semester, multi-course dataset including Business Calculus, Linear Algebra, and Calculus III. By comparing unproctored homework and proctored exam performance, we analyze how student learning behaviors shift in AI-accessible environments, particularly through epistemic off-loading of mathematical work. Guided by a sociocognitive framework, we employ complementary measures -- performance gaps, homework-exam correlations, and Wasserstein distance -- to characterize divergence between practice and mastery. Results reveal a growing integrity gap as course content shifts from procedural to conceptual and spatially intensive mathematics. In both Business Calculus and Linear Algebra, differences in homework format (online versus hand-written, TA-graded) do not yield substantively different performance patterns, indicating that paper-based homework is not inherently more resistant to AI-mediated offloading. While homework retains partial predictive validity in procedural courses, upper-division courses exhibit a collapse in alignment between homework and exams, indicating that unproctored assessments increasingly reflect synthetic fluency rather than internalized understanding. These findings highlight the need to rethink assessment practices in the AI era.}




\maketitle

\section{Introduction}\label{sec1}
The rapid development of large language models (LLMs) such as ChatGPT, Claude, and Gemini, alongside mature symbolic computation systems, has fundamentally altered the undergraduate mathematics classroom. Tasks that once defined the daily labor of students, such as algebraic manipulation, matrix computation, differentiation and integration, and even routine multivariable calculus operations, can now be executed instantly, accurately, and with minimal user input. As these tools become widely accessible, their presence is no longer peripheral but structurally embedded in how students engage with mathematics both inside and outside the classroom.

This technological shift has generated a growing tension in mathematics education. On one hand, artificial intelligence can function as a scaffolding tool \cite{wood1976role}, supporting learning by offering step-by-step explanations, visualizations, and immediate feedback. Used judiciously, such tools may lower cognitive barriers, reduce anxiety, and allow students to focus on higher-level reasoning. On the other hand, AI can operate as a prosthetic, completing core mathematical work on behalf of students and potentially bypassing the processes through which understanding is traditionally constructed. Distinguishing between these two roles, AI as aid versus AI as substitute, has become increasingly difficult, particularly in unsupervised learning environments.

This tension strikes at the heart of longstanding debates in collegiate mathematics about instructional goals. Traditional curricula have balanced procedural fluency with conceptual understanding, and computation with modeling, abstraction, and structural reasoning. Generative and symbolic AI systems now perform many procedural tasks at a level comparable to or exceeding that of undergraduates \cite{frieder2023mathematical}, raising a fundamental question: if machines can reliably execute the procedures we teach, what should count as mathematical knowledge and competence for students?


Existing research on generative AI in higher education has primarily examined tool adoption, student attitudes, academic integrity concerns, or short-term performance effects \cite{sullivan2023chatgpt,fuchs2023exploring,xie2023ai}. While these studies provide important insights into the risks of AI-assisted cheating and the limitations of detection-based countermeasures, they largely frame AI as an external disruption to be managed rather than as a force that fundamentally reshapes the epistemic function of coursework and assessment. Recent conceptual work has begun to challenge this view. Li and Wang \cite{haifengli2023} argue that in the era of generative AI, traditional homework has lost its educational meaning, as task completion can be fully delegated to AI systems, necessitating a redesign of assignments and evaluation toward human -- AI collaboration and process-oriented assessment. Similarly, Leaton Gray et al. \cite{leaton2025ai} contend that AI does not merely facilitate misconduct but exposes deeper structural weaknesses in contemporary assessment regimes, calling for a shift from surveillance-driven integrity models toward validity-centered, ethically grounded pedagogy. However, empirical evidence quantifying how AI availability reshapes the alignment between formative and summative assessment across different types of mathematics courses remains limited.

This study addresses that gap through a two-semester, multi-course empirical investigation spanning Business Calculus, Linear Algebra, and Calculus III. These courses were deliberately chosen to represent a progression from procedurally dominated mathematics to conceptually structured and spatially intensive content. Using detailed assessment data, we examine how student performance evolves in the presence of readily available AI tools, with particular attention to differences between unproctored assessments (homework) and proctored assessments (exams). This comparison allows us to explore whether an integrity gap emerges or widens over time, reflecting shifts in student learning behaviors rather than isolated incidents of misuse.

The study is guided by the following research questions:
\begin{itemize}[
  label={},
  leftmargin=3em,
  labelwidth=2.5em,
  labelsep=0.5em,
  align=left
]
\item[\textbf{RQ1}] How does the availability of AI tools affect student performance across different types of undergraduate mathematics courses?

\item[\textbf{RQ2}] Do these effects differ between procedurally focused courses and conceptually or spatially structured courses?

\item[\textbf{RQ3}] How do discrepancies between homework and exam performance reflect shifts in student learning strategies and reliance on AI?

\item[\textbf{RQ4}] What implications do these patterns have for curriculum design, assessment practices, and the future goals of undergraduate mathematics education?

\end{itemize}

By positioning AI not as a pedagogical gadget but as a disruptive epistemic force, this work aims to contribute empirical evidence to a rapidly evolving conversation. Understanding how AI reshapes what students do, what they learn, and what assessments actually measure is essential for rethinking mathematics education in the AI era.

\section{Theoretical Framework}\label{sec2}
To rigorously analyze the impact of generative AI on student learning and assessment outcomes, this study adopts a sociocognitive lens that treats AI not merely as a neutral technological aid, but as an epistemic actor whose educational meaning depends on how it is appropriated by students. We synthesize the theory of instrumental genesis, the distinction between pragmatic and epistemic values, and the didactical contract to explain how AI reshapes mathematical activity.
\subsection{The Artifact versus The Instrument: Genesis and ZPD}
We begin by distinguishing between AI systems (e.g., ChatGPT) as artifacts and as instruments. An artifact is a material or symbolic tool made available to the learner; it does not, by itself, guarantee learning. According to the theory of instrumental genesis \cite{verillon1995cognition}, an artifact becomes an instrument for learning only through a developmental process in which the learner integrates the tool into their cognitive activity.


Instrumental genesis involves two intertwined processes. Instrumentation refers to the ways in which the artifact shapes the user's thinking and behavior. For example, when AI automatically provides coordinate-based solutions, symbolic manipulations, or completed derivations, students may cease attempts to visualize vectors, reason geometrically, or anticipate qualitative behavior. In this sense, the tool reorganizes the student's cognitive schema by redefining what constitutes a legitimate mathematical action.

Instrumentalization, by contrast, refers to how learners adapt and shape the tool to serve mathematical goals, learning how to constrain prompts, request explanations aligned with formal definitions, or use AI outputs as objects for verification and critique. Instrumentalization is essential if AI is to function as a learning instrument rather than a solution generator.

Our central claim is that for many students in this study, instrumental genesis remained incomplete. The AI did not become a mathematical instrument but instead functioned as a prosthetic artifact, performing the mathematical work on the student's behalf \cite{risko2016cognitive}. In these cases, instrumentation dominated instrumentalization: the AI replaced core cognitive labor rather than scaffolding it. As a result, task completion, particularly in unproctored settings, cannot be assumed to reflect internalized mathematical understanding.

This phenomenon is closely related to Skemp's distinction \cite{skemp1976relational} between instrumental understanding (rule-following without reasons) and relational understanding (knowing both what to do and why). Generative AI enables what may be termed \textit{synthetic fluency}: students can produce correct outputs and polished solutions without engaging in the reasoning processes traditionally required to develop relational understanding. This synthetic fluency threatens core mathematical practices, including reasoning, representation, and abstraction.

Traditionally, the Zone of Proximal Development (ZPD) \cite{vygotsky1978mind} represents the gap between independent ability and ability with guidance. We posit that AI may artificially inflate the ZPD. By instantly solving problems beyond the student's actual capability, AI mimics competence, creating an integrity gap. The student appears capable within the ZPD during unproctored homework, but collapses to their actual independent level during proctored exams, revealing that the scaffolding was actually a prosthetic.



\subsection{Pragmatic versus Epistemic Values}
To explain why this effect varies by course, we adopt Artigue's distinction between the \textit{pragmatic} and \textit{epistemic} value of tools \cite{artigue2002learning}.
Pragmatic value refers to a tool's effectiveness in helping students get the job done, typically by producing correct answers efficiently. Epistemic value refers to a tool's capacity to support conceptual understanding, meaning-making, and the construction of mathematical knowledge.

This plays out differently across the curriculum: In Business Calculus (high pragmatic utility), Tasks here are often algorithmic (e.g., solving systems of equations). AI offers high pragmatic value by efficiently handling the drudgery of calculation. While this risks bypassing procedural skill acquisition, the cognitive demand is often low enough that students can still verify the output, potentially using AI as an efficiency tool.

Linear Algebra (mixed regime) occupies a mixed epistemic regime characterized by high levels of abstraction and the need to coordinate symbolic, geometric, and structural representations. While AI systems can efficiently perform matrix computations and routine algebraic manipulations, they provide limited support for the development of relational understanding required for reasoning about vector spaces, linear transformations, and proofs. In this setting, AI often functions as a translation or execution tool rather than a conceptual scaffold, enabling correct symbolic output without guaranteeing comprehension of the underlying structures. As a result, student reliance on AI in Linear Algebra is highly unstable: it can remain benign for computation-oriented tasks but becomes epistemically disruptive when assessments shift toward conceptual or structural reasoning.

Calculus III (the epistemic black box) presents a unique barrier: the combination of complex algebra with high-order 3D spatial reasoning. Contrary to the assumption that visual tasks require human intervention, the sheer difficulty of the material often exceeds the student's independent cognitive capacity. Consequently, this course is highly prone to black box off-loading. When students cannot visualize a flux integral or parameterize a surface, they do not partial-scaffold; they surrender the entire problem to AI. Here, reliance is driven not by convenience, but by necessity, creating the highest risk of the prosthetic effect where epistemic value is eliminated entirely.

This distinction leads directly to the study's central hypothesis: the largest gaps between unproctored homework and proctored exams will occur where the pragmatic value of AI is high but epistemic engagement is low. Homework performance may reflect successful pragmatic use of AI, while exam performance reveals the absence of internalized conceptual schemes.

\subsection{Disruption of the Didactical Contract}
Finally, we situate these dynamics within Brousseau's Didactical Contract \cite{brousseau2002theory}. Traditionally, assigning a problem carries the implicit understanding that solving it requires engaging with the underlying mathematics.

Generative AI disrupts this contract by satisfying the visible condition of the task-producing a correct answer -- without satisfying its educational purpose. This aligns with Brousseau's Topaze Effect, where the answer is obtained through means that circumvent the intended learning process. In the presence of AI, this effect is systematized: the student provides the correct output, but the mathematical knowledge remains external, residing entirely within the artifact.

This framework suggests that AI functions as an epistemic actor that redefines human-essential work. If students offload both the procedural manipulation and the conceptual modeling to the AI, they abdicate their role in the epistemic process, rendering assessment scores a measure of the tool's capability rather than the student's understanding.

As a result, assessment artifacts, especially homework, may no longer function as reliable indicators of learning. Divergence between homework and exam performance in this study is therefore interpreted not merely as an issue of academic integrity, but as evidence of a deeper rupture in the didactical contract, in which traditional signals of competence no longer correspond to actual cognitive development.

\section{Methodology}\label{sec3}

This study was conducted over a two-semester period across three undergraduate mathematics courses chosen to represent distinct positions along a cognitive and epistemic spectrum: Business Calculus, Linear Algebra, and Calculus III. These courses were selected deliberately to examine whether the impact of AI-mediated problem solving varies systematically with the nature of mathematical knowledge and practices emphasized.


These courses form a progression from procedurally dominated mathematics to conceptually and representationally demanding mathematics, providing a natural testbed for examining differential patterns of AI reliance and learning outcomes.

\subsection{Participants}

The student population across all three courses was ethnically homogeneous, consisting entirely of Chinese students, which reduces cross-cultural variability in prior educational background. However, substantial within-population variation exists due to China's differentiated high school mathematics tracks. In the Chinese secondary education system, students are streamed into science-track and liberal-arts-track programs, which involve markedly different levels of mathematical intensity and rigor. 

Students entering Business Calculus, particularly those majoring in Accounting and Management, predominantly come from the liberal-arts track and typically have weaker formal preparation in advanced algebra and calculus, whereas students enrolled in Linear Algebra and Calculus III are more likely to come from the science track with substantially stronger pre-university mathematics training. Importantly, these background differences are stable across semesters and therefore cannot account for the observed collapse in homework-exam alignment in upper-division courses.

\begin{itemize}
\item Business Calculus students were predominantly first-year undergraduates majoring in Accounting, Finance, Management, and related business disciplines. Class sizes were typically around 50 students per section.

\item Linear Algebra students consisted primarily of first-year Data Science majors and second-year Computer Science majors, with a small number of Accounting and Finance majors and a few senior students retaking the course. Typical class sizes ranged from 25 to 30 students.


\item Calculus III students were primarily second-year Data Science majors, again with a small number of senior retake students. Class sizes were similar to Linear Algebra, generally between 25 and 30 students.
\end{itemize}

The same instructor taught the same courses included in the study, ensuring consistency in instructional style, assessment design, and grading standards.

\subsection{Data Collection}






Assessment data was collected through two primary modalities, reflecting standard undergraduate assessment practices.

\begin{itemize}
\item Automated Assessment (WebWorK) \cite{roth2008evaluating,toth2013measuring,prat2021webwork}: All three courses utilized WebWorK, an online homework delivery system. This platform generates randomized numerical parameters for each problem, ensuring that while the method remains consistent across students, the specific answer varies. This prevents simple answer-sharing but does not prevent the use of AI solvers that can process the specific parameters. Calculus III relied solely on WebWorK for homework assessment.

\item Manual Assessment (Hand-Graded): Business Calculus and Linear Algebra included a hand-graded homework component. This required students to submit written workings, theoretically demanding a demonstration of the process. However, modern AI tools (e.g., GPT-5.2 with vision capabilities) can now generate step-by-step written solutions, making this modality vulnerable to transcription-based academic dishonesty.

\end{itemize}

In all courses, midterm and final examinations were conducted in-person under proctored conditions, without access to AI tools or external computational aids beyond those explicitly permitted (e.g., basic calculators when allowed).

\subsection{Variables}

The primary independent distinction in this study is between unproctored and proctored assessment environments. Homework (online and written) is treated as unproctored, with high potential for AI use. While collaboration and external resources were not prohibited, AI usage was neither explicitly encouraged nor systematically restricted. Midterm and final exams are treated as proctored, with low or negligible AI access. These assessments are intended to measure students' independent mathematical competence under time constraints. The primary dependent measures are students' normalized scores on homework assignments and exams within each course.

\subsection{Measures}




Guided by the theoretical framework in Section \ref{sec2}, we pre-specified a set of complementary analytical measures designed to capture different ways in which the relationship between practice and mastery may break down in the presence of generative AI. Rather than relying on a single summary statistic, we adopt a multi-dimensional divergence framework that was defined prior to data analysis and applied uniformly across courses and semesters.

This approach is motivated by two considerations. First, the theory of instrumental genesis predicts that AI may differentially affect task completion and conceptual mastery. Second, the pragmatic-epistemic distinction suggests that disruptions may manifest not only as average score differences, but also as weakened alignment and structural fragmentation between formative and summative assessments. Accordingly, three outcome measures were defined a priori.

\subsubsection{Performance Delta}

The first measure is the raw performance delta, defined for each student as:
$$\Delta = \text{Average Score}_{\text{Homework}} - \text{Score}_{\text{Exam}}$$
This metric captures the most direct form of divergence between unproctored and proctored performance. It reflects the extent to which students' success in AI-accessible environments translates into independent performance under controlled conditions.

We emphasize that $\Delta$ is treated as a baseline descriptor rather than a definitive indicator of learning loss. It was included a priori to provide continuity with existing assessment comparisons and to serve as a reference point for the additional measures described below.

\subsubsection{Predictive Validity of Homework}

The second measure evaluates the predictive validity of homework by computing the Spearman rank correlation coefficient ($r$) between students' homework averages and exam scores within each course and semester. Spearman correlations were used to reduce sensitivity to ceiling effects in homework scores.

In conventional instructional settings, homework is expected to function as guided practice, resulting in a positive association between homework and exam performance. Drawing on this assumption, we use the strength of this correlation as an indicator of whether homework remains aligned with independent mastery.

A reduction in this association -- referred to here as correlation collapse -- is interpreted as evidence that homework performance no longer reliably signals conceptual understanding. Importantly, this measure was selected in advance to detect decoupling effects predicted by prosthetic AI use, independent of overall score levels or distributional shape.

\subsubsection{Distributional Morphology}

To capture structural differences between unproctored and proctored performance beyond central tendency, we pre-specified a distributional morphology measure based on the Wasserstein distance. This metric quantifies how much “mass” must be transported to transform one score distribution into another, making it sensitive to shifts in spread, skewness, and polarization.

Let $F_{H}$ and $F_{E}$ denote the empirical Cumulative Distribution Functions (CDFs) of homework and exam scores, respectively, within a given course and semester.
Given the discrete nature of our grading scale ($k \in \{0, 1, \dots, 100\}$), we compute the Wasserstein distance \cite{villani2008optimal} $W_1$ using the $L_1$ norm of the CDFs:
\begin{align*}
W_1(H, E) = \sum_{k=0}^{100} \left| F_{H}(k) - F_{E}(k) \right|.
\end{align*}
The Wasserstein distance $W_1$ captures the cost of transforming the homework distribution into the exam distribution. This metric is particularly sensitive to polarization, allowing us to distinguish between a general increase in difficulty (a uniform shift) and a fracture in the student body (bimodality), where one group retains mastery and another collapses.

This measure was selected a priori for three reasons. First, unlike mean differences, the Wasserstein distance is robust to ceiling effects in homework scores, which are expected under high pragmatic AI use. Second, it is sensitive to heterogeneous changes, including increased dispersion or the emergence of subpopulations, without requiring assumptions about unimodality or normality. Third, it provides a single, interpretable scalar that summarizes distributional divergence while preserving information about the full score distribution.

In this framework, a small Wasserstein distance indicates that homework and exam performances are distributionally aligned, consistent with homework functioning as effective preparation for independent mastery. A large distance indicates substantial misalignment between practice and assessment, consistent with AI-mediated inflation of unproctored performance and fragmentation of mastery under proctored conditions.

Importantly, the Wasserstein distance is used as a complementary diagnostic alongside the raw performance delta and homework-exam correlation. While $\Delta$ captures magnitude and correlation captures predictive alignment, the Wasserstein distance captures structural divergence between assessment regimes. 


By defining these three measures prior to analysis and grounding each in a specific theoretical expectation, this study minimizes the risk of post hoc metric selection. The measures are not treated as interchangeable alternatives, but as orthogonal diagnostics that jointly characterize magnitude, alignment, and structure in student performance. All subsequent analyses apply this framework consistently across courses and semesters.

\section{Data Analysis \& Results}\label{sec4}

In this section, we present the empirical results across the three studied courses. For each course, we analyze the divergence between unproctored (homework) and proctored (exam) performance using the three metrics defined in our methodology: the raw performance delta ($\Delta$), the predictive validity correlation ($r$), and the distributional morphology ($W_1$). All data analysis and generation of figures were performed in R.

\subsection{Business Calculus}

The data for Business Calculus were collected from the semester of Spring 2025. The results for Business Calculus characterize a course where AI usage appears to be predominantly pragmatic -- used for efficiency rather than total cognitive replacement.

\begin{figure}[htbp]
  \centering
  \includegraphics[width=0.8\textwidth]{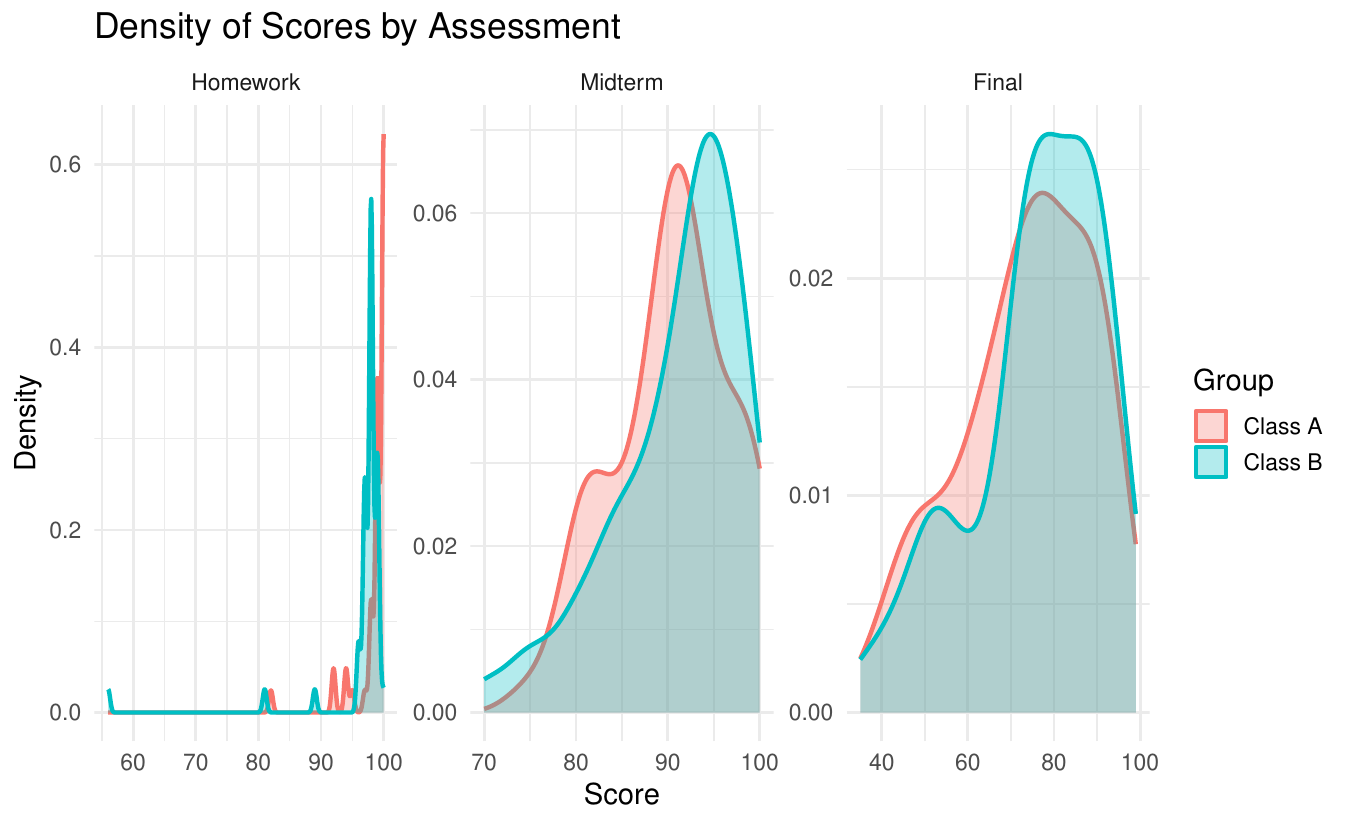}
  \caption{Density plots for two groups of Business Calculus}
  \label{fig:density2400}
\end{figure}

Figure~\ref{fig:density2400} shows the distributions of assessment scores for Business Calculus. Homework in Class A was delivered via the WebWorK platform, while homework in Class B was hand-graded by a teaching assistant. Although both classes received identical problem sets, the WebWorK assignments included randomized numerical parameters to mitigate answer sharing.

\begin{table}[htbp]
\centering
\caption{Mean $\pm$ SD for Business Calculus}
\label{tab:mean-sd-2400}
\begin{tabular}{lccc}
\toprule
Group &
Homework &
Midterm &
Final \\
\midrule
Class A &
$98.5 \pm 3.0$ &
$90.4 \pm 6.3$ &
$73.4 \pm 15.3$ \\
Class B &
$96.6 \pm 6.5$ &
$91.1 \pm 7.0$ &
$75.7 \pm 15.4$ \\
\bottomrule
\end{tabular}
\end{table}

As shown in Table \ref{tab:mean-sd-2400}, students in both sections achieved near-perfect averages on homework ($98.5$ for Class A; $96.6$ for Class B). SD in the table denotes standard deviation. However, performance dropped significantly on the proctored final exam ($73.4$ and $75.7$, respectively). Class A has $\Delta=25.1$ points and Class B has $\Delta=20.9$ points.
While this gap is substantial, it reflects a substantial but relatively uniform decline typical of courses where students use tools to bypass tedious calculation but retain procedural awareness. Midterm exams, on the other hand, admit relative small $\Delta$'s.

\begin{table}[htbp]
\centering
\caption{Spearman Correlations for Business Calculus}
\label{tab:correlations2400}
\begin{tabular}{l
                S[table-format=1.3]
                S[table-format=1.3]
                S[table-format=1.3]}
\toprule
Group   & {Homework--Midterm } & {Homework--Final } & {Midterm--Final} \\
\midrule
Class A   &  0.508   &  0.410  &   0.604 \\
Class B   &  0.484   &  0.546  &   0.495 \\
\bottomrule
\end{tabular}
\end{table}

Despite the drop in raw scores, the relationship between homework and exam performance remained statistically significant. Table \ref{tab:correlations2400} reports Spearman correlations between homework and final exam scores of $r = 0.410$ (Class A) and $r = 0.546$ (Class B). This moderate positive correlation suggests that students who performed better on homework generally performed better on exams, implying that homework exercises still retained some value as predictive practice.

\begin{table}[htbp]
\centering
\caption{Wasserstein Distance for Business Calculus}
\label{tab:wasserstein2400}
\sisetup{round-mode=places, round-precision=1}
\begin{tabular}{l
                S[table-format=2.1]
                S[table-format=2.1]
                S[table-format=2.1]}
\toprule
Group   & {Homework--Midterm} & {Homework--Final} & {Midterm--Final} \\
\midrule
Class A & 8.08 & 25.2 & 17.1 \\
Class B & 7.04 & 20.9 & 16.1 \\
\bottomrule
\end{tabular}
\end{table}

The Wasserstein distance ($W_1$) between homework and final exam distributions was calculated at 25.2 (Class A) and 20.9 (Class B). As seen in the density plots (Figure \ref{fig:density2400}), the exam distribution shifted leftward (lower scores) but retained a relatively unimodal shape. This uniform shift suggests that while the material was harder without tools, the class performed somewhat homogeneously, without a deep fracture in the student population.

An additional observation from the Business Calculus results is that the mode of homework administration-online versus paper-based does not produce a substantively different performance pattern. Despite differences in format, with Class A completing homework via WebWorK and Class B submitting hand-written, TA-graded solutions, both groups exhibited similarly high homework averages, comparable homework-exam correlations, and nearly identical Wasserstein distances between homework and exam score distributions. 

This suggests that, in a procedurally oriented course such as Business Calculus, the assessment medium itself does not significantly mitigate or exacerbate AI-mediated performance inflation. Even when students are required to submit written work, contemporary AI tools can generate step-by-step solutions that are easily transcribed, rendering paper-based homework no more resistant to epistemic offloading than online platforms. Consequently, the observed integrity gap in Business Calculus appears to be driven primarily by the cognitive nature of the tasks rather than by the homework delivery format.

\subsection{Linear Algebra}

The data for Linear Algebra were collected from the semesters of Spring 2025 and Fall 2025. The results for Linear Algebra represent a mixed regime, effectively bridging the gap between the pragmatic stability of Business Calculus and the prosthetic collapse of Calculus III. This course, which demands both computational fluency (matrix operations) and abstract reasoning (vector spaces, proofs), exhibits highly variable alignment between unproctored and proctored performance. Linear Algebra emerges as the most diagnostically revealing course in the dataset, exhibiting extreme sensitivity to both cohort characteristics and the epistemic nature of assessment tasks.

Linear Algebra serves as the critical transition point -- or stress test -- of this study's central hypothesis. Unlike Business Calculus, which is predominantly procedural, and Calculus III, which is heavily conceptual and spatial, Linear Algebra occupies an intermediate epistemic position. Its assessment outcomes fluctuate sharply depending on whether evaluative emphasis is placed on computational fluency or on conceptual and structural reasoning.

This split is empirically visible across semesters. In Spring 2025, when the final exam emphasized computation-oriented tasks, performance divergence between unproctored homework and proctored exams remained comparable to Business Calculus, with integrity gaps $\Delta$ typically in the range of 20--36 points. By contrast, in Fall 2025, when the final exam shifted toward concept-oriented problems, the integrity gap widened dramatically. For Classes E, F, and G, $\Delta$ exceeded 48 points, matching or surpassing the collapse observed in Calculus III.

This sharp discontinuity demonstrates that Linear Algebra is not merely intermediate in content, but epistemically unstable: small shifts away from procedural computation toward conceptual abstraction are sufficient to trigger a breakdown in homework-exam alignment. As such, Linear Algebra provides the clearest evidence that AI-mediated integrity gaps are driven not by course level alone, but by the epistemic nature of the assessed tasks.

\begin{figure}[htbp]
  \centering
  \includegraphics[width=0.8\textwidth]{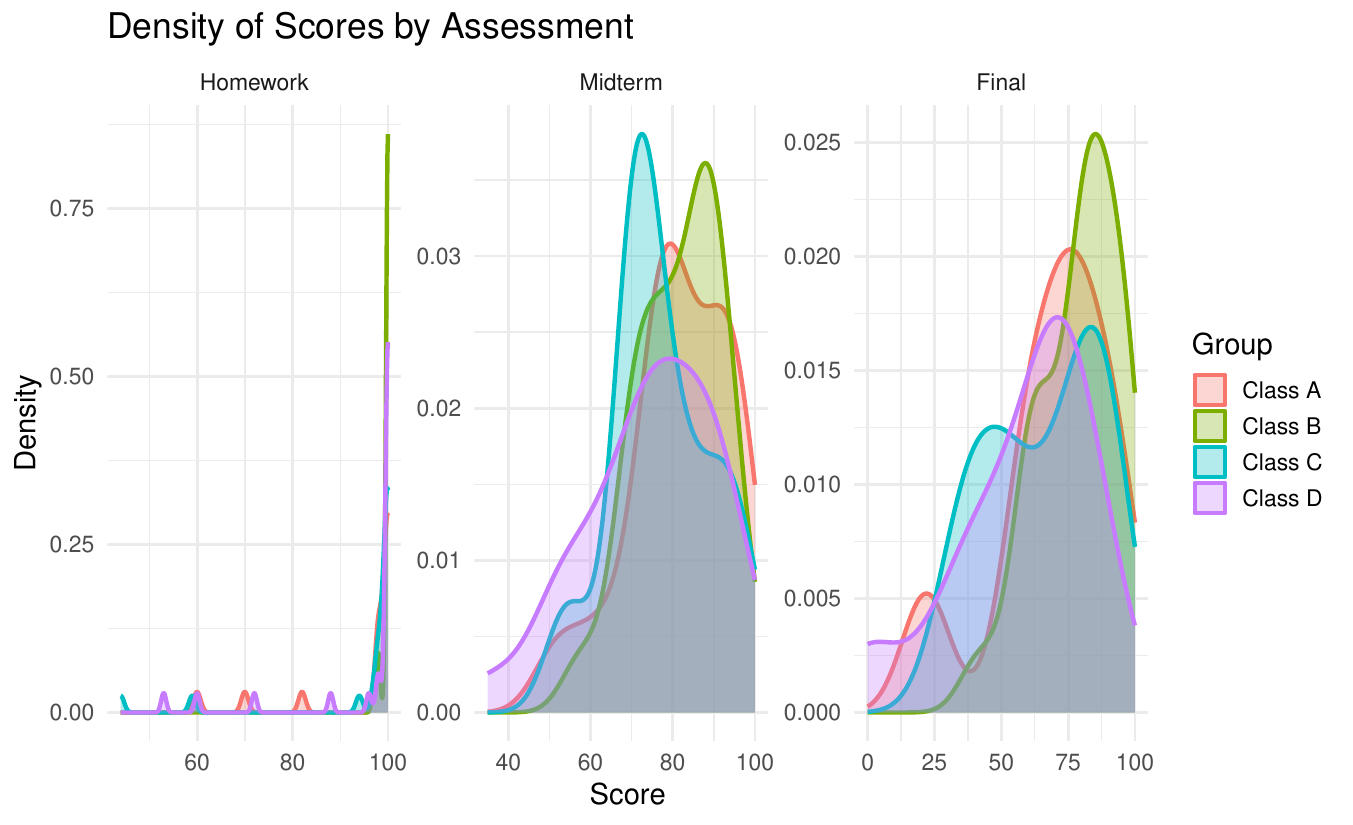}
  \caption{Density plots for four groups of Linear Algebra (2025 Spring)}
  \label{fig:density2995spring}
\end{figure}

\begin{figure}[htbp]
  \centering
  \includegraphics[width=0.8\textwidth]{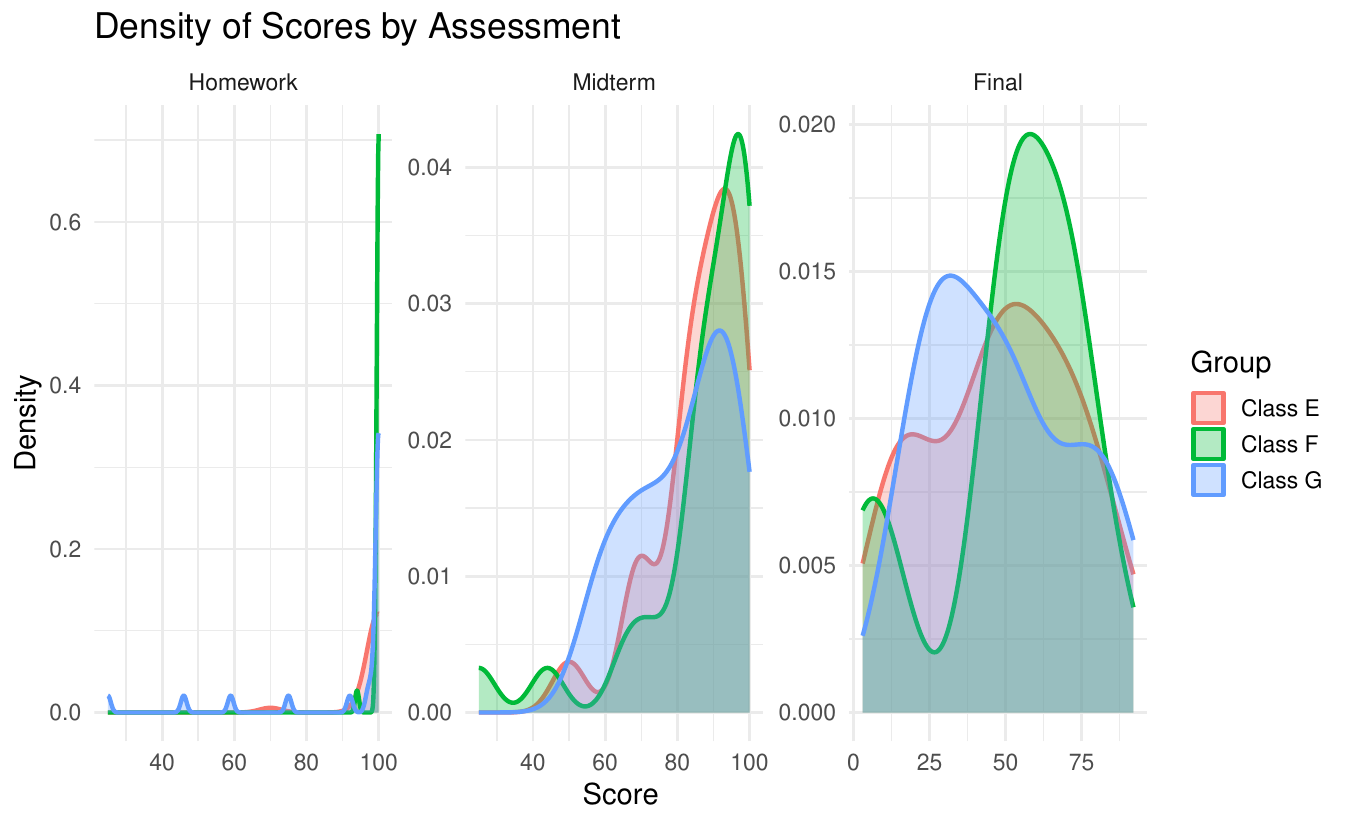}
  \caption{Density plots for three groups of Linear Algebra (2025 Fall)}
  \label{fig:density2995fall}
\end{figure}

Figure~\ref{fig:density2995spring} and Figure~\ref{fig:density2995fall} show the distributions of assessment scores for Linear Algebra from Spring 2025 and Fall 2025, respectively. Homework in Classes B and D were delivered via the WebWorK platform, while homework in Classes A and C were hand-graded by the instructor. The problem sets for Classes A and C were drawn from textbook exercises, whereas the WebWorK assignments included randomized numerical parameters. In Fall 2025, homework in all classes (Classes E, F, and G) were delivered via the WebWorK platform.


The homework assignments and midterm problems are consistent within each semester, and student performance on homework and midterm exams is relatively comparable. However, the final exam in Spring 2025 is computation-oriented and the final exam in Fall 2025 is concept-oriented, which leads to that the final grades in Fall 2025 are substantially lower than those in Spring 2025 (Table \ref{tab:mean-sd-2995}).

\begin{table}[htbp]
\centering
\caption{Mean $\pm$ SD for Linear Algebra by Semester}
\label{tab:mean-sd-2995}
\begin{tabular}{llccc}
\toprule
Semester & Group & Homework & Midterm & Final \\
\midrule
\multirow{4}{*}{2025 Spring}
 & Class A & $94.3 \pm 12.0$ & $81.6 \pm 12.5$ & $69.6 \pm 21.9$ \\
 & Class B & $99.7 \pm 0.80$ & $81.7 \pm 10.0$ & $79.2 \pm 15.1$ \\
 & Class C & $94.8 \pm 14.3$ & $77.1 \pm 12.1$ & $65.7 \pm 21.0$ \\
 & Class D & $94.8 \pm 12.7$ & $74.6 \pm 15.8$ & $58.2 \pm 25.0$ \\
\midrule
\multirow{3}{*}{2025 Fall}
 & Class E & $98.3 \pm 6.30$  & $86.3 \pm 12.3$ & $48.7 \pm 24.5$ \\
 & Class F & $99.8 \pm 1.15$ & $86.7 \pm 17.9$ & $51.2 \pm 24.7$ \\
 & Class G & $91.3 \pm 19.9$ & $81.5 \pm 13.7$ & $48.8 \pm 23.6$ \\
\bottomrule
\end{tabular}
\end{table}

Similar to the other courses, homework averages remained near the ceiling across all groups, ranging from 91.3 to 99.8 (Table \ref{tab:mean-sd-2995}). However, the integrity gap widened significantly compared to Business Calculus. In the Spring 2025 cohort (Classes A, B, C, and D), final exam (computation-oriented) averages ranged from 58.2 to 79.2, resulting in performance deltas ($\Delta$) of approximately 20--36 points, which is similar to the case of Business Calculus. In the Fall 2025 cohort (Classes E, F, and G), this gap became catastrophic: while homework scores remained high (98--99), final exam (concept-oriented) averages plummeted to 48--51, creating deltas exceeding 48 points, which is similar to the case of Calculus III. This stark contrast suggests that while some cohorts managed to bridge the gap, others experienced a prosthetic failure comparable to the upper-division Calculus III students.

\begin{table}[htbp]
\centering
\caption{Spearman Correlations for Linear Algebra by Semester}
\label{tab:correlation2995}
\begin{tabular}{ll
                S[table-format=1.3]
                S[table-format=1.3]
                S[table-format=1.3]}
\toprule
Semester & Group 
& {Homework--Midterm} 
& {Homework--Final} 
& {Midterm--Final} \\
\midrule
\multirow{4}{*}{2025 Spring}
 & Class A &  0.718  &  0.461  &  0.649 \\
 & Class B &  0.209  &  0.372  &  0.553 \\
 & Class C & -0.061  & -0.176  &  0.722 \\
 & Class D &  0.581  &  0.516  &  0.613 \\
\midrule
\multirow{3}{*}{2025 Fall}
 & Class E &  0.641  &  0.471  &  0.577 \\
 & Class F &  0.305  &  0.252  &  0.420 \\
 & Class G &  0.582  &  0.283  &  0.256 \\
\bottomrule
\end{tabular}
\end{table}

The predictive power of homework in Linear Algebra is best described as unstable. Unlike Business Calculus, where correlations were consistently moderate, Linear Algebra displays extreme variability. 
From Spring 2025 (Table \ref{tab:correlation2995} Classes A, B, C, and D, computation-oriented final exam), Class A maintained a moderate homework-final correlation ($r=0.461$), but Class C exhibited a negative correlation ($r=-0.176$), implying that students with higher homework scores actually performed worse on independent assessments -- a hallmark of the Topaze Effect where reliance on tools actively hinders learning. 
From Fall 2025 (Table \ref{tab:correlation2995} Classes E, F, and G, concept-oriented final exam), the correlations weakened further, with Class F and G showing low associations ($r \approx$ 0.25--0.28), signaling a breakdown in the assessment's validity.

\begin{table}[htbp]
\centering
\caption{Wasserstein Distance for Linear Algebra by Semester}
\label{tab:wasserstein2995}
\sisetup{round-mode=places, round-precision=1}
\begin{tabular}{ll
                S[table-format=2.1]
                S[table-format=2.1]
                S[table-format=2.1]}
\toprule
Semester & Group
& {Homework--Midterm}
& {Homework--Final}
& {Midterm--Final} \\
\midrule
\multirow{4}{*}{2025 Spring}
 & Class A & 12.6 & 24.6 & 12.0 \\
 & Class B & 18.0 & 20.5 &  4.6 \\
 & Class C & 18.5 & 29.0 & 11.7 \\
 & Class D & 20.2 & 36.6 & 16.3 \\
\midrule
\multirow{3}{*}{2025 Fall}
 & Class E & 12.1 & 49.7 & 37.6 \\
 & Class F & 13.1 & 48.0 & 34.9 \\
 & Class G & 13.5 & 42.5 & 32.7 \\
\bottomrule
\end{tabular}
\end{table}

The Wasserstein distance metrics confirm this transitional status. In Spring 2025 (Table \ref{tab:wasserstein2995} Classes A, B, C, and D), $W_1$ values ranged from 20.5 to 36.6, generally higher than Business Calculus but lower than Calculus III, more resembling Business Calculus. However, the Fall 2025 data (Table \ref{tab:wasserstein2995} Classes E, F, and G) shows $W_1$ values spiking to 42.5--49.7. This scale of values is comparable to the scale in Calculus III. As visually confirmed by the density plots (Figure \ref{fig:density2995spring} is in the pattern of Figure \ref{fig:density2400}  \& Figure \ref{fig:density2995fall} is more like Figure \ref{fig:density3415}), the homework distributions remain sharp peaks at 100, while the exam distributions flatten and spread into the failing range. This indicates that Linear Algebra sits on a fault line: depending on the cohort and specific content (e.g., computational vs. abstract), the course can tip from a manageable tool-assisted environment into a state of total epistemic offloading.

Sample problems of the final exam for 2025 Spring and Fall are included in the Appendix~\ref{secA1}.

\subsection{Calculus III}
The data for Calculus III were collected from the semester of Fall 2025 and the homework is administered on WebWorK only. 
In contrast to Business Calculus, the results for Calculus III reveal a fundamental breakdown in the relationship between unproctored and proctored assessment, consistent with the Prosthetic Artifact hypothesis.

\begin{figure}[htbp]
  \centering
  \includegraphics[width=0.8\textwidth]{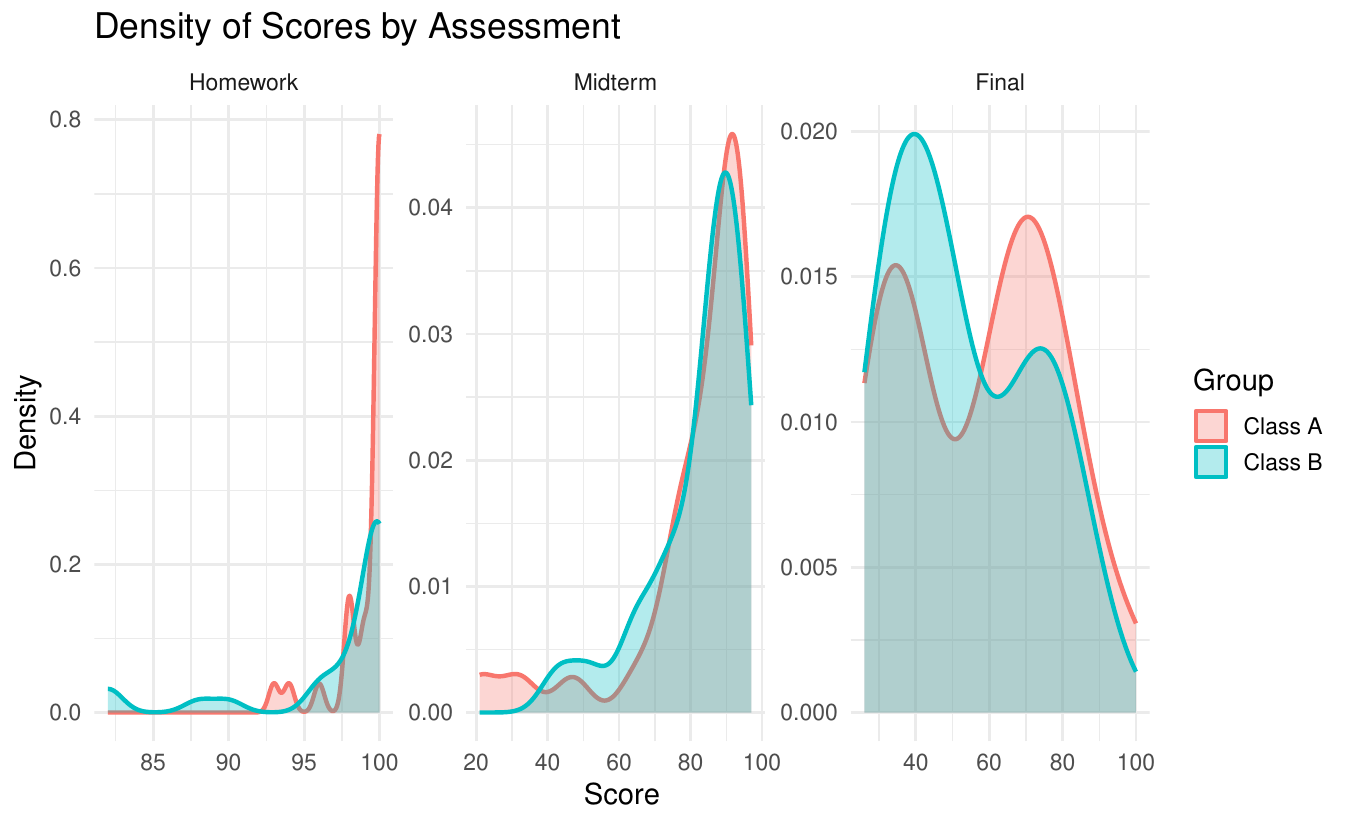}
  \caption{Density plots for two groups of Calculus III}
  \label{fig:density3415}
\end{figure}

Figure~\ref{fig:density3415} demonstrates a pronounced bimodal and polarized distribution for both Class A and Class B in Calculus III. A subset of students maintains moderate exam performance, suggesting retained independent competence, whereas a large fraction of the cohort performs poorly under proctored conditions. This pattern is consistent with heterogeneous modes of AI engagement, ranging from constrained support to extensive cognitive offloading.

\begin{table}[htbp]
\centering
\caption{Mean $\pm$ SD for Calculus III}
\label{tab:mean-sd-3415}
\begin{tabular}{lccc}
\toprule
Group &
Homework &
Midterm &
Final \\
\midrule
Class A &
$99.1 \pm 1.8$ &
$81.6 \pm 18.4$ &
$57.3 \pm 21.3$ \\
Class B &
$97.1 \pm 5.4$ &
$82.5 \pm 13.7$ &
$52.7 \pm 19.2$ \\
\bottomrule
\end{tabular}
\end{table}

Table \ref{tab:mean-sd-3415} shows that homework averages remained at the ceiling ($99.1$ for Class A; $97.1$ for Class B), statistically indistinguishable from the lower-division course. However, the proctored final exam scores collapsed to failing averages. Class A has an average of $57.3$ with $\Delta = 41.8$ points. Class B has an average of $52.7$ with $\Delta = 44.4$ points. This divergence is nearly double that observed in Business Calculus, indicating that for spatially complex material, the removal of the AI tool resulted in a severe loss of capability.

\begin{table}[htbp]
\centering
\caption{Spearman Correlations for Calculus III}
\label{tab:correlation3415}
\begin{tabular}{l
                S[table-format=1.3]
                S[table-format=1.3]
                S[table-format=1.3]}
\toprule
Group   & {Homework--Midterm } & {Homework--Final } & {Midterm--Final} \\
\midrule
Class A & 0.290 & 0.273 & 0.693 \\
Class B & 0.188 & 0.391 & 0.535 \\
\bottomrule
\end{tabular}
\end{table}

The predictive validity of homework effectively vanished in this course. Table \ref{tab:correlation3415} reports homework-final correlations of $r = 0.273$ (Class A) and $r = 0.391$ (Class B). In Class A, the correlation between homework and midterm was as low as $0.290$. This weak association implies that a student's homework grade is no longer a reliable indicator of their independent mastery of the material.

\begin{table}[htbp]
\centering
\caption{Wasserstein Distance for Calculus III}
\label{tab:wasserstein3415}
\sisetup{round-mode=places, round-precision=1}
\begin{tabular}{l
                S[table-format=2.1]
                S[table-format=2.1]
                S[table-format=2.1]}
\toprule
Group   & {Homework--Midterm} & {Homework--Final} & {Midterm--Final} \\
\midrule
Class A  &   17.5  &   41.7  &   24.9 \\
Class B  &   14.6  &   44.4  &   29.8 \\
\bottomrule
\end{tabular}
\end{table}

The structural divergence is most visible in the Wasserstein distance, which spiked to 41.7 (Class A) and 44.4 (Class B) -- among the highest values observed in the study. The density plots in Figure \ref{fig:density3415} visually confirm this fracture. While the homework distribution (far left panel) is a sharp peak at 100, the final exam distribution (far right panel) flattens completely, spreading across the entire range with significant density in the failing region ($<60$). This massive transport of ``probability mass'' confirms that the high homework scores were artificial, masking a severe deficit in the conceptual understanding required for the proctored environment. 

\section{Discussion}\label{sec5}
The results of this study provide empirical support for the hypothesis that the impact of AI on student learning is not uniform but is mediated by the epistemic value of the course content. By contrasting the pragmatic environment of Business Calculus with the epistemic demands of Calculus III, we observe two distinct patterns of instrumental genesis.
\subsection{The Prosthetic Effect and the Integrity Gap}
This subsection addresses RQ1 and RQ2 by examining how AI-related performance divergence manifests differently across procedurally oriented and conceptually demanding mathematics courses.
The most immediate finding is the magnitude of the performance divergence between unproctored and proctored assessments. In Business Calculus, we observed a significant but manageable drop from homework averages ($97$) to final exam averages ($74$). 
While this indicates reliance on external aids, the moderate correlations between homework and exams ($r \approx 0.50$) suggest that homework still functioned, at least partially, as a form of practice.
Here, AI likely served a pragmatic role -- automating calculation while leaving some procedural structure intact.

In the middle of this spectrum lies Linear Algebra, which functions as a transitional domain. The data reveals an unstable landscape where the integrity gap varies wildly. In some sections (Spring Class B), the gap mirrored the manageable levels of Business Calculus ($\Delta \approx 20$). In others (Fall Classes E, F, G), the gap widened to match or even exceed the catastrophic levels of Calculus III ($\Delta > 48$). This variability supports the hypothesis that Linear Algebra contains both pragmatic tasks (which AI handles efficiently without destroying learning) and abstract conceptual tasks (where AI acts as a prosthetic). When the course leans too heavily into the latter without proctoring, the scaffolding collapses.

In contrast, Calculus III revealed a severe integrity gap. Despite maintaining near-perfect homework averages ($98$), students' independent performance collapsed on the Final Exam ($55$).
This 40-point divergence suggests that for spatially and conceptually complex tasks, AI functioned not as a scaffold but as a prosthetic artifact. 
The tool supplanted core cognitive processes rather than supporting their development. This confirms our hypothesis that courses requiring high-order visualization and conceptual modeling are paradoxically more vulnerable to total offloading than rote procedural courses.

An important cross-course finding concerns the limited role of homework format in mitigating AI-mediated performance divergence. In both Business Calculus and Linear Algebra, sections using online homework platforms (WebWorK) and those requiring hand-written, hand-graded submissions exhibited comparable homework averages, similar homework-exam correlations, and nearly indistinguishable distributional divergence measures. This pattern suggests that paper-based homework, often assumed to enforce greater epistemic engagement, does not meaningfully restore assessment validity in the presence of contemporary AI tools. In Linear Algebra in particular, where hand-written solutions are traditionally viewed as essential for conceptual learning, the integrity gap persisted despite the requirement to submit full written work. AI-generated step-by-step solutions can be readily transcribed, rendering the assessment medium largely orthogonal to the underlying issue of cognitive offloading. Consequently, the observed breakdown in homework-exam alignment is better explained by the epistemic demands of the task than by the mode of homework administration.

\subsection{The Erosion of Take-Home Assessment Validity}
This subsection addresses RQ3 by analyzing how discrepancies between homework and exam performance reflect shifts in student engagement, learning strategies, and reliance on AI.
Perhaps the most alarming finding is the correlation collapse observed in the upper-division data. In traditional pedagogy, homework performance is a predictor of exam readiness. Our data shows this relationship has severed in the AI era for conceptual courses. 

In Business Calculus, the correlation between homework and midterm scores remained statistically relevant ($r \approx 0.50$). In Calculus III, this predictive validity vanished ($r \approx$ 0.19--0.29). This lack of correlation implies that in Calculus III, a student's homework score is no longer a measure of their own capability, but rather a measure of the tool's capability. Linear Algebra provides the ``smoking gun'' for this decay in validity. The observation of a negative correlation ($r = -0.176$) in Class C is particularly revealing. In this instance, high homework performance was not merely uninformative; it was a counter-signal to mastery. This suggests a dependency trap, where the most diligent users of AI for homework were the least equipped for independent reasoning. Unlike Business Calculus, where the contract held, and Calculus III, where it severed, Linear Algebra shows the contract fraying -- validity is no longer guaranteed and fluctuates unpredictably.

This represents a fundamental rupture in the Didactical Contract. The students fulfilled the clause of the contract (submitting correct answers) without fulfilling the spirit (engaging in the epistemic process), resulting in a Topaze Effect where correct answers mask a complete absence of learning.

\subsection{Distributional Fracture}
From a distributional perspective, this analysis further supports RQ3 by showing that AI-mediated learning loss is structural rather than uniform.
The analysis of distributional morphology via the Wasserstein distance further illuminates the structural nature of this learning loss. The distance between homework and final exam distributions in Calculus III ($W_1 \approx$ 42-44) was nearly double that of Business Calculus ($W_1 \approx$ 21--25). 

The distributional morphology of Linear Algebra further illustrates this transition. The Wasserstein distances ($W_1$) in the Fall 2025 cohort (42--49) actually surpassed those of Calculus III, indicating extreme polarization. This suggests that even in mixed courses, the student population is fracturing into two distinct groups: those who use the tool to learn (retaining the distribution shape) and those who use the tool to replace labor (flattening the distribution).

This high Wasserstein distance in Calculus III reflects a polarization of the student body. The density plots reveal that while homework scores remained clustered at the ceiling (synthetic fluency), exam scores flattened and shifted drastically downward. This indicates that the middle ground of partial understanding is disappearing; students either master the material independently or surrender the cognitive labor entirely to the artifact, leaving them helpless during proctored assessments.

\subsection{Implications for Curriculum Design} 
This subsection directly addresses RQ4 by translating observed assessment misalignments into implications for curriculum design and assessment practices in the AI era.
The data suggests that AI-proofing assignments through difficulty alone is a failed strategy. In fact, increasing the conceptual difficulty (as in Linear Algebra and Calculus III) appeared to increase the reliance on prosthetic AI use. For Procedural Courses (Business Calculus), the focus should shift from calculation to interpretation. Since AI excels at the algorithmic execution, assessments should verify if students can contextualize the AI's output. For Conceptual Courses (Linear Algebra and Calculus III), the traditional take-home problem set is effectively dead as a summative metric. If unproctored work is retained, it should be graded on process (e.g., video explanations, oral defenses) rather than product.

\section{Conclusion}\label{sec6}
This study set out to trace the impact of generative AI on undergraduate mathematics education across a cognitive spectrum. Our findings indicate that AI has ceased to be a peripheral tool and has become a central epistemic actor that fundamentally alters the nature of student engagement.


Across courses (RQ1), we observed a systematic divergence between unproctored homework performance and proctored exam performance, indicating that AI accessibility fundamentally alters what assessment scores represent. However, these effects were not uniform. The magnitude, structure, and educational meaning of this divergence varied markedly by course type (RQ2). In procedurally oriented courses such as Business Calculus, AI use appeared largely pragmatic: homework retained moderate predictive validity for exams, suggesting partial alignment between practice and mastery. In contrast, in conceptually and spatially demanding courses such as Calculus III, AI use became predominantly prosthetic, producing near-perfect homework scores alongside severe exam collapse.


These discrepancies between homework and exam performance (RQ3) reveal a deeper breakdown in student engagement and learning strategies. The collapse of homework-exam correlations and the large Wasserstein distances observed in upper-division courses indicate that unproctored assessments increasingly measure synthetic fluency rather than internalized understanding. In such settings, homework performance reflects the capability of the tool rather than the competence of the student, signaling a rupture in the didactical contract.

Taken together, these findings carry significant implications for curriculum design and assessment in the AI era (RQ4). Increasing task difficulty or reverting to traditional paper-based homework alone does not mitigate AI reliance and may, in some cases, intensify epistemic offloading. Across both Business Calculus and Linear Algebra, hand-written homework failed to preserve alignment between practice and independent mastery, demonstrating that assessment validity cannot be restored through format-level changes alone. To preserve the integrity of mathematical assessment, unproctored homework should be reclassified as formative practice rather than summative evidence of mastery, while greater emphasis should be placed on proctored and process-oriented assessments that require explanation, reasoning, and conceptual accountability.

These results challenge the assumption that AI serves merely as a ``tutor'' or ``scaffold'' within the ZPD. Instead, for many students, AI artificially inflates the ZPD, creating an illusion of competence that vanishes immediately upon the removal of the tool. To preserve the integrity of mathematics education, we recommend a structural pivot in assessment strategies:

\begin{itemize}
\item De-weighting Unproctored Tasks: Homework should be treated as formative practice (graded for completion or low-stakes) rather than a measure of mastery.

\item Elevation of Proctored Environments: Summative assessments must return to controlled environments to ensure the measurement of human cognition rather than synthetic fluency.

\item Process-Oriented Assessment: In courses where conceptual understanding is paramount, students must be assessed on their ability to explain why a result is true, a task that currently remains harder to offload completely than the generation of the result itself.
\end{itemize}

Ultimately, AI is not inherently harmful to mathematics education, but it is disruptive. As our data on the integrity gap demonstrates, continuing to rely on pre-AI assessment models in a post-AI world does not measure student learning -- it increasingly measures the capability of the machine rather than the competence of the learner.

\section{Limitations and Future Work}
Several limitations of this study should be acknowledged. First, the data were drawn from a single institution with a relatively homogeneous student population, which limits the generalizability of the findings across institutional types, and cultural contexts. Second, this study is observational rather than experimental; while the divergence between unproctored and proctored assessments is theoretically grounded and empirically robust, causal claims about individual students' AI usage cannot be made without direct measures of tool interaction. In particular, the absence of direct AI-usage logs prevents fine-grained attribution of individual performance patterns to specific modes of AI engagement. Third, homework scores are used as a proxy for AI-accessible performance, which captures aggregate effects but does not distinguish between productive, scaffolded use of AI and complete cognitive offloading at the individual level. 

Future work should therefore incorporate mixed-method designs, including controlled experiments, AI-usage logging, think-aloud protocols, and student interviews, to disentangle how different modes of AI engagement shape learning. Longitudinal studies spanning additional courses and institutions are also needed to determine whether the observed integrity gap stabilizes, widens, or attenuates as both students and curricula adapt to AI-mediated mathematics learning.

\bmhead{Acknowledgements}
This work is partially supported by Wenzhou Kean University and Department of Education of Zhejiang Province under contract JGBA2024540 and JGBA2024546.

\section*{Declarations}

\noindent\textbf{Competing Interests} \\
The authors declare no conflicts of interest (financial or non-financial).

\vspace{0.5em}
\noindent\textbf{Funding} \\
This work is partially supported by Wenzhou Kean University and Department of Education of Zhejiang Province under contract JGBA2024540 and JGBA2024546.

\vspace{0.5em}
\noindent\textbf{Ethics Approval and Consent to Participate} \\
This study did not involve medical or psychological experiments on human participants and did not require formal ethics approval. All data analyzed were generated as part of routine educational practice and were anonymized prior to analysis.

\vspace{0.5em}
\noindent\textbf{Consent for Publication} \\
Not applicable.

\vspace{0.5em}
\noindent\textbf{Data Availability} \\
The datasets generated and analyzed during the current study are available from the corresponding author upon reasonable request.

\vspace{0.5em}
\noindent\textbf{Code Availability} \\
The code used for data analysis is available from the corresponding author upon reasonable request.

\vspace{0.5em}
\noindent\textbf{Materials Availability} \\
Not applicable.

\vspace{0.5em}
\noindent\textbf{Author Contributions} \\
All authors contributed to the conception and design of the study. Data collection and analysis were performed by the authors. All authors drafted and critically revised the manuscript and approved the final version for submission.








\begin{appendices}

\section{Samples of final exams of Linear Algebra}\label{secA1}


\subsection{Spring 2025 (computation-oriented)}

\begin{enumerate}
     \item Suppose that $A = \begin{bmatrix} 3 &2&-1\\x&-2&2\\3&y&-1 \end{bmatrix}$ has an eigenvector $\textbf{v} = \begin{bmatrix} 1\\-2\\3 \end{bmatrix}$. Find the values of $x$, $y$ and the eigenvalue corresponding to $\textbf{v}$. 
     \item Suppose that matrix $C$ is similar to the diagonal matrix $B= \begin{bmatrix} 1 & & \\ & -1 & \\ & & &1 \end{bmatrix}$. Find $C^2$. 
 \end{enumerate}

\subsection{Fall 2025 (concept-oriented)}

Let $A$ be a $3 \times 3$ matrix, and suppose $\textbf{v}_1$, $\textbf{v}_2$, $\textbf{v}_3$ are eigenvectors of $A$ corresponding to eigenvalues 0, 1, 2. Let $\textbf{u} = \textbf{v}_1 + \textbf{v}_2 + \textbf{v}_3$.
     \begin{enumerate}
         \item Determine whether the set $\{\textbf{u}, A\textbf{u}, A^2\textbf{u}\}$ is linearly independent or linearly dependent. Justify your answer.
         \item Find the values of $a$, $b$, and $c$ such that $A$ is similar to $B = \begin{bmatrix} 0&0&a \\ 1&0&b \\ 0&1&c \end{bmatrix}$. 
     \end{enumerate}

\end{appendices}


\bibliography{sn-bibliography}

\begin{thebibliography}{}
\renewcommand{\doi}[1]{\url{https://doi.org/#1}}
\bibcommenthead

\bibitem [\protect \citeauthoryear {%
Artigue%
}{%
Artigue%
}{%
{\protect \APACyear {2002}}%
}]{%
artigue2002learning}
\APACinsertmetastar {%
artigue2002learning}%
\begin{APACrefauthors}%
Artigue, M.%
\end{APACrefauthors}%
\unskip\
\newblock
\APACrefYearMonthDay{2002}{}{}.
\newblock
{\BBOQ}\APACrefatitle {Learning mathematics in a CAS environment: The genesis
  of a reflection about instrumentation and the dialectics between technical
  and conceptual work} {Learning mathematics in a cas environment: The genesis
  of a reflection about instrumentation and the dialectics between technical
  and conceptual work}.{\BBCQ}
\newblock
\APACjournalVolNumPages{International journal of computers for mathematical
  learning}{7}{3}{245--274,}
\newblock

\newblock

\PrintBackRefs{\CurrentBib}

\bibitem [\protect \citeauthoryear {%
Brousseau%
}{%
Brousseau%
}{%
{\protect \APACyear {2002}}%
}]{%
brousseau2002theory}
\APACinsertmetastar {%
brousseau2002theory}%
\begin{APACrefauthors}%
Brousseau, G.%
\end{APACrefauthors}%
\unskip\
\newblock
\APACrefYear{2002}.
\newblock
\APACrefbtitle {Theory of didactical situations in mathematics: Didactique des
  math{\'e}matiques, 1970--1990} {Theory of didactical situations in
  mathematics: Didactique des math{\'e}matiques, 1970--1990}.
\newblock
\APACaddressPublisher{Dordrecht, The Netherlands}{Springer}.
\PrintBackRefs{\CurrentBib}

\bibitem [\protect \citeauthoryear {%
Frieder%
\ \protect \BOthers {.}}{%
Frieder%
\ \protect \BOthers {.}}{%
{\protect \APACyear {2023}}%
}]{%
frieder2023mathematical}
\APACinsertmetastar {%
frieder2023mathematical}%
\begin{APACrefauthors}%
Frieder, S.%
, Pinchetti, L.%
, Chevalier, A.%
, Griffiths, R\BHBI R.%
, Salvatori, T.%
, Lukasiewicz, T.%
\BDBL {}Berner, J.%
\end{APACrefauthors}%
\unskip\
\newblock
\APACrefYearMonthDay{2023}{}{}.
\newblock
{\BBOQ}\APACrefatitle {Mathematical Capabilities of Chat{GPT}} {Mathematical
  capabilities of chat{GPT}}.{\BBCQ}
\newblock
 \APACrefbtitle {Thirty-seventh Conference on Neural Information Processing
  Systems Datasets and Benchmarks Track.} {Thirty-seventh conference on neural
  information processing systems datasets and benchmarks track.}
\PrintBackRefs{\CurrentBib}

\bibitem [\protect \citeauthoryear {%
Fuchs%
}{%
Fuchs%
}{%
{\protect \APACyear {2023}}%
}]{%
fuchs2023exploring}
\APACinsertmetastar {%
fuchs2023exploring}%
\begin{APACrefauthors}%
Fuchs, K.%
\end{APACrefauthors}%
\unskip\
\newblock
\APACrefYearMonthDay{2023}{}{}.
\newblock
{\BBOQ}\APACrefatitle {Exploring the opportunities and challenges of NLP models
  in higher education: is Chat GPT a blessing or a curse?} {Exploring the
  opportunities and challenges of nlp models in higher education: is chat gpt a
  blessing or a curse?}{\BBCQ}
\newblock
 \APACrefbtitle {Frontiers in education} {Frontiers in education}\ (\BVOL~8,
  \BPG~1166682).
\PrintBackRefs{\CurrentBib}

\bibitem [\protect \citeauthoryear {%
Leaton~Gray%
, Edsall%
\BCBL {}\ \BBA {} Parapadakis%
}{%
Leaton~Gray%
\ \protect \BOthers {.}}{%
{\protect \APACyear {2025}}%
}]{%
leaton2025ai}
\APACinsertmetastar {%
leaton2025ai}%
\begin{APACrefauthors}%
Leaton~Gray, S.%
, Edsall, D.%
\BCBL {} Parapadakis, D.%
\end{APACrefauthors}%
\unskip\
\newblock
\APACrefYearMonthDay{2025}{}{}.
\newblock
{\BBOQ}\APACrefatitle {AI-Based Digital Cheating At University, and the Case
  for New Ethical Pedagogies} {Ai-based digital cheating at university, and the
  case for new ethical pedagogies}.{\BBCQ}
\newblock
\APACjournalVolNumPages{Journal of Academic Ethics}{}{}{1--18,}
\newblock

\newblock

\PrintBackRefs{\CurrentBib}

\bibitem [\protect \citeauthoryear {%
Li%
\ \BBA {} Wang%
}{%
Li%
\ \BBA {} Wang%
}{%
{\protect \APACyear {2023}}%
}]{%
haifengli2023}
\APACinsertmetastar {%
haifengli2023}%
\begin{APACrefauthors}%
Li, H.%
\BCBT {}\ \BBA {} Wang, W.%
\end{APACrefauthors}%
\unskip\
\newblock
\APACrefYearMonthDay{2023}{}{}.
\newblock
{\BBOQ}\APACrefatitle {Design and Evaluation of Student Assignment in the Era
  of Generative Artificial Intelligence (Chinese)} {Design and evaluation of
  student assignment in the era of generative artificial intelligence
  (chinese)}.{\BBCQ}
\newblock
\APACjournalVolNumPages{Open Education Research}{29}{3}{,}
\newblock

\newblock

\PrintBackRefs{\CurrentBib}

\bibitem [\protect \citeauthoryear {%
Prat%
\ \BBA {} Code%
}{%
Prat%
\ \BBA {} Code%
}{%
{\protect \APACyear {2021}}%
}]{%
prat2021webwork}
\APACinsertmetastar {%
prat2021webwork}%
\begin{APACrefauthors}%
Prat, A.%
\BCBT {}\ \BBA {} Code, W.J.%
\end{APACrefauthors}%
\unskip\
\newblock
\APACrefYearMonthDay{2021}{}{}.
\newblock
{\BBOQ}\APACrefatitle {WeBWorK log files as a rich source of data on student
  homework behaviours} {Webwork log files as a rich source of data on student
  homework behaviours}.{\BBCQ}
\newblock
\APACjournalVolNumPages{International Journal of Mathematical Education in
  Science and Technology}{52}{10}{1540--1556,}
\newblock

\newblock

\PrintBackRefs{\CurrentBib}

\bibitem [\protect \citeauthoryear {%
Risko%
\ \BBA {} Gilbert%
}{%
Risko%
\ \BBA {} Gilbert%
}{%
{\protect \APACyear {2016}}%
}]{%
risko2016cognitive}
\APACinsertmetastar {%
risko2016cognitive}%
\begin{APACrefauthors}%
Risko, E.F.%
\BCBT {}\ \BBA {} Gilbert, S.J.%
\end{APACrefauthors}%
\unskip\
\newblock
\APACrefYearMonthDay{2016}{}{}.
\newblock
{\BBOQ}\APACrefatitle {Cognitive offloading} {Cognitive offloading}.{\BBCQ}
\newblock
\APACjournalVolNumPages{Trends in cognitive sciences}{20}{9}{676--688,}
\newblock

\newblock

\PrintBackRefs{\CurrentBib}

\bibitem [\protect \citeauthoryear {%
Roth%
, Ivanchenko%
\BCBL {}\ \BBA {} Record%
}{%
Roth%
\ \protect \BOthers {.}}{%
{\protect \APACyear {2008}}%
}]{%
roth2008evaluating}
\APACinsertmetastar {%
roth2008evaluating}%
\begin{APACrefauthors}%
Roth, V.%
, Ivanchenko, V.%
\BCBL {} Record, N.%
\end{APACrefauthors}%
\unskip\
\newblock
\APACrefYearMonthDay{2008}{}{}.
\newblock
{\BBOQ}\APACrefatitle {Evaluating student response to WeBWorK, a web-based
  homework delivery and grading system} {Evaluating student response to
  webwork, a web-based homework delivery and grading system}.{\BBCQ}
\newblock
\APACjournalVolNumPages{Computers \& Education}{50}{4}{1462--1482,}
\newblock

\newblock

\PrintBackRefs{\CurrentBib}

\bibitem [\protect \citeauthoryear {%
Skemp%
}{%
Skemp%
}{%
{\protect \APACyear {1976}}%
}]{%
skemp1976relational}
\APACinsertmetastar {%
skemp1976relational}%
\begin{APACrefauthors}%
Skemp, R.R.%
\end{APACrefauthors}%
\unskip\
\newblock
\APACrefYearMonthDay{1976}{}{}.
\newblock
{\BBOQ}\APACrefatitle {Relational understanding and instrumental understanding}
  {Relational understanding and instrumental understanding}.{\BBCQ}
\newblock
\APACjournalVolNumPages{Mathematics teaching}{77}{1}{20--26,}
\newblock

\newblock

\PrintBackRefs{\CurrentBib}

\bibitem [\protect \citeauthoryear {%
Sullivan%
, Kelly%
\BCBL {}\ \BBA {} McLaughlan%
}{%
Sullivan%
\ \protect \BOthers {.}}{%
{\protect \APACyear {2023}}%
}]{%
sullivan2023chatgpt}
\APACinsertmetastar {%
sullivan2023chatgpt}%
\begin{APACrefauthors}%
Sullivan, M.%
, Kelly, A.%
\BCBL {} McLaughlan, P.%
\end{APACrefauthors}%
\unskip\
\newblock
\APACrefYearMonthDay{2023}{}{}.
\newblock
{\BBOQ}\APACrefatitle {ChatGPT in higher education: Considerations for academic
  integrity and student learning} {Chatgpt in higher education: Considerations
  for academic integrity and student learning}.{\BBCQ}
\newblock
\APACjournalVolNumPages{Journal of Applied Learning \& Teaching}{6}{1}{31--40,}
\newblock

\newblock

\PrintBackRefs{\CurrentBib}

\bibitem [\protect \citeauthoryear {%
Toth%
}{%
Toth%
}{%
{\protect \APACyear {2013}}%
}]{%
toth2013measuring}
\APACinsertmetastar {%
toth2013measuring}%
\begin{APACrefauthors}%
Toth, P.F.%
\end{APACrefauthors}%
\unskip\
\newblock
\APACrefYearMonthDay{2013}{}{}.
\newblock
{\BBOQ}\APACrefatitle {Measuring efficiency of teaching mathematics online:
  experiences with WeBWorK} {Measuring efficiency of teaching mathematics
  online: experiences with webwork}.{\BBCQ}
\newblock
\APACjournalVolNumPages{Procedia-Social and Behavioral
  Sciences}{89}{}{276--282,}
\newblock

\newblock

\PrintBackRefs{\CurrentBib}

\bibitem [\protect \citeauthoryear {%
Verillon%
\ \BBA {} Rabardel%
}{%
Verillon%
\ \BBA {} Rabardel%
}{%
{\protect \APACyear {1995}}%
}]{%
verillon1995cognition}
\APACinsertmetastar {%
verillon1995cognition}%
\begin{APACrefauthors}%
Verillon, P.%
\BCBT {}\ \BBA {} Rabardel, P.%
\end{APACrefauthors}%
\unskip\
\newblock
\APACrefYearMonthDay{1995}{}{}.
\newblock
{\BBOQ}\APACrefatitle {Cognition and artifacts: A contribution to the study of
  though in relation to instrumented activity} {Cognition and artifacts: A
  contribution to the study of though in relation to instrumented
  activity}.{\BBCQ}
\newblock
\APACjournalVolNumPages{European journal of psychology of
  education}{}{}{77--101,}
\newblock

\newblock

\PrintBackRefs{\CurrentBib}

\bibitem [\protect \citeauthoryear {%
Villani%
}{%
Villani%
}{%
{\protect \APACyear {2008}}%
}]{%
villani2008optimal}
\APACinsertmetastar {%
villani2008optimal}%
\begin{APACrefauthors}%
Villani, C.%
\end{APACrefauthors}%
\unskip\
\newblock
\APACrefYear{2008}.
\newblock
\APACrefbtitle {Optimal transport: old and new} {Optimal transport: old and
  new}\ (\BVOL~338).
\newblock
\APACaddressPublisher{Berlin, Heidelberg}{Springer}.
\PrintBackRefs{\CurrentBib}

\bibitem [\protect \citeauthoryear {%
Vygotsky%
}{%
Vygotsky%
}{%
{\protect \APACyear {1978}}%
}]{%
vygotsky1978mind}
\APACinsertmetastar {%
vygotsky1978mind}%
\begin{APACrefauthors}%
Vygotsky, L.S.%
\end{APACrefauthors}%
\unskip\
\newblock
\APACrefYear{1978}.
\newblock
\APACrefbtitle {Mind in society: The development of higher psychological
  processes} {Mind in society: The development of higher psychological
  processes}\ (\BVOL~86).
\newblock
\APACaddressPublisher{Cambridge, MA}{Harvard university press}.
\PrintBackRefs{\CurrentBib}

\bibitem [\protect \citeauthoryear {%
Wood%
, Bruner%
\BCBL {}\ \BBA {} Ross%
}{%
Wood%
\ \protect \BOthers {.}}{%
{\protect \APACyear {1976}}%
}]{%
wood1976role}
\APACinsertmetastar {%
wood1976role}%
\begin{APACrefauthors}%
Wood, D.%
, Bruner, J.S.%
\BCBL {} Ross, G.%
\end{APACrefauthors}%
\unskip\
\newblock
\APACrefYearMonthDay{1976}{}{}.
\newblock
{\BBOQ}\APACrefatitle {The role of tutoring in problem solving} {The role of
  tutoring in problem solving}.{\BBCQ}
\newblock
\APACjournalVolNumPages{Journal of child psychology and
  psychiatry}{17}{2}{89--100,}
\newblock

\newblock

\PrintBackRefs{\CurrentBib}

\bibitem [\protect \citeauthoryear {%
Xie%
, Wu%
\BCBL {}\ \BBA {} Chakravarty%
}{%
Xie%
\ \protect \BOthers {.}}{%
{\protect \APACyear {2023}}%
}]{%
xie2023ai}
\APACinsertmetastar {%
xie2023ai}%
\begin{APACrefauthors}%
Xie, Y.%
, Wu, S.%
\BCBL {} Chakravarty, S.%
\end{APACrefauthors}%
\unskip\
\newblock
\APACrefYearMonthDay{2023}{}{}.
\newblock
{\BBOQ}\APACrefatitle {AI meets AI: Artificial intelligence and academic
  integrity-A survey on mitigating AI-assisted cheating in computing education}
  {Ai meets ai: Artificial intelligence and academic integrity-a survey on
  mitigating ai-assisted cheating in computing education}.{\BBCQ}
\newblock
 \APACrefbtitle {Proceedings of the 24th annual conference on information
  technology education} {Proceedings of the 24th annual conference on
  information technology education}\ (\BPGS\ 79--83).
\PrintBackRefs{\CurrentBib}

\end{thebibliography}

\end{document}